\documentclass{svjour3}                     % onecolumn (standard format)
\smartqed  % flush right qed marks, e.g. at end of proof
\usepackage{latexsym}
\usepackage{amssymb}
\usepackage{amsmath}
\usepackage{rotating}
\usepackage[misc]{ifsym}
\newtheorem{them}{Theorem}[section]

\newtheorem{pro}{Proposition}[section]

\journalname{JOTA}
\begin{document}

\title{Nonsmooth Steepest Descent Method by Proximal Subdifferentials in Hilbert Spaces}
%\thanks{Communicated by Boris S. Mordukhovich}
\titlerunning{Nonsmooth steepest descent method by proximal subdifferentials in Hilbert spaces}



\author{ Zhou Wei \and Qing Hai He}
        %etc.


\institute{Zhou Wei(\Letter) \at Department of Mathematics, Yunnan University,
 Kunming 650091, P. R. China\\ \email{wzhou@ynu.edu.cn}
  \and Qing Hai He \at Department of Mathematics, Yunnan University,
 Kunming 650091, P. R. China\\ \email{heqh@ynu.edu.cn}}
                    %  \\

\date{Received: date / Accepted: date}

\maketitle

\begin{abstract}
In this paper, we first study nonsmooth steepest descent method for nonsmooth functions defined on Hilbert space and establish the corresponding algorithm by proximal subgradients. Then, we use this algorithm to find stationary points for those functions satisfying prox-regularity and Lipschitzian continuity. As one application, the established algorithm is used to search the minimizer of lower semicontinuous and convex functions on finite-dimensional space. The  convergent theorem, as one extension and improvement of the existing converging result for twice continuously differentiable convex functions, is also presented therein.


\end{abstract}

\section{Introduction}
Given a smooth function on an Euclidean space, one kind of optimization problem is to find the global/local minimizers of the given function. This kind of problems has been receiving a great deal of attention and it is of significance to solve it by using the iteration of line search methods. For solving this problem, it is used generally to find the steepest descent unit direction, and then the given function can be reduced along this direction and finally attain the global minimizer hopefully. This so-called steepest descent method nowadays has been applied to various application areas such as numerical optimization, computing mathematics, probabilistic analysis, algorithm, differential equations, vector optimization and so on (cf. \cite{3,4,CY,DS,FS,17,P,23} and references therein). This method could be traced back to Cauchy and subsequently studied and improved by many authors. One of most useful improvement is the conjugate gradient method. It is known that Hestenes and Stiefel \cite{10} studied the linear conjugate gradient method for solving linear systems with positive definite coefficient matrices in 1950s, and afterwards in 1960s, the first nonlinear conjugate gradient method, as one of the earliest known techniques for solving large-scale nonlinear optimization problems, was introduced by Fletcher and Reeves in \cite{9}. The classic algorithm for smooth functions by steepest descent method is presented in Section 2.

Note that the main idea of the steepest descent method is, by using the steepest descent unit direction, to produce a gradient sequence that converges to zero with mild assumptions such as Zoutendijk, Wolfe or Goldstein conditions. It cannot be guaranteed that this method  converges to a minimizer, but only that it is attracted by stationary points. Naturally, it is interesting and important to establish the corresponding algorithm for nonsmooth functions on Hilbert space along the line by this method. Motivated by it, in this paper, we mainly study the nonsmooth version of steepest descent method by proximal subdifferential of nonsmooth functions defined on Hilbert space, and aim to construct one appropriate algorithm being used to find the stationary point by the nonsmooth steepest descent method.

The paper is organized as follows. In Section 2, we will give some definitions and preliminary through this paper. Our notation is basically standard and conventional in the area of variational analysis. Section 3 is devoted to establishing the algorithm for nonsmooth functions on Hilbert space by nonsmooth steepest descent method. This algorithm is used to find stationary points of those functions with prox-regularity and Lipschitzian continuity assumptions in Hilbert space. The newly-established algorithm is also applied to the search of the minimizer for lower semicontinuous  G\^{a}teaux differentiable convex functions on the finite-dimensional space (see Theorem 3.4). The conclusion of this paper is presented in Section 4.

\section{Preliminaries} Let $H$ be a {\em Hilbert space} equipped with the {\em inner product}
$\langle\cdot, \cdot\rangle$ and the corresponding {\em norm} $\|\cdot\|$. Given a Hilbert space $H$ and a multifunction $F: H\rightrightarrows H$, the symbol
\begin{equation*}
\begin{array}r
\mathop{\rm Limsup}\limits_{y\rightarrow x}F(x):=\Big\{\zeta\in H: \exists \ {\rm sequences} \ x_n\rightarrow x \ {\rm and} \ \zeta_n\stackrel{w}\longrightarrow \zeta \ {\rm with}\  \\
\zeta_n\in F(x_n) \ {\rm for\ all \ } n\in \mathbb{N} \Big\}
\end{array}
\end{equation*}
signifies the {\it sequential Painlev\'{e}-Kuratowski outer/upper limit} of $F(x)$ as $y\rightarrow x$.

For $x\in H$ and $\delta>0$, let $B(x, \delta)$ denote the {\em open ball} with center $x$ and radius $\delta$. Given a set $A\subset H$, denote $\overline A$ and ${\rm int}(A)$ the {\it norm closure} and the {\it interior } of $A$, respectively. If a sequence $\{v_n\}$ converges weakly to $v$, it is denoted by $v_n\stackrel{w}\longrightarrow v$.

For any point $z\in H$, the distance between $z$ and $S$ is given by
$$
d(z, S):=\inf\{\|z-s\|\,:\,s\in S\}.
$$
The set of {\em closest points} to $z$ in $S$ is denoted by
\begin{equation*}
P_S(z):=\{s\in S\,:\,\|z-s\|=d(z, S)\}.
\end{equation*}
Let $x\in S$. Recall that the {\it proximal normal cone} of $S$ at $x$, denoted by $N^p(S, x)$, is defined as:
$$
N^p(S, x):=\big\{\zeta\in H: \zeta=t(z-x)\  {\rm for \ some}\  (z, t)\in H\times (0, +\infty) \ {\rm with} \ x\in P_S(z)\big\}.
$$
It is known from \cite{8} that $\zeta\in N^p(S, x)$ if and only if there exist two constants $\sigma,\delta\in ]0, +\infty[$ such that
\begin{equation*}
\langle\zeta, s-x\rangle\leq \sigma\|s-x\|^2,\,\,\,\forall s\in
B(x, \delta)\cap S.
\end{equation*}

The {\it limiting(Mordukhovich) normal cone} of $S$ at $x$, denoted by $N(S, x)$, is defined as
\begin{equation*}
N(S, x):=\mathop{\rm Limsup}_{y\stackrel{S}\longrightarrow x}N^p(S, y),
\end{equation*}
where
$y\stackrel{S}\longrightarrow x$ means $y\rightarrow x$ and $y\in S$. Thus, $\zeta\in N(S, x)$ if and only if there exists a sequence $\{(x_n, \zeta_n)\}$ in $S\times H$ such that $x_n\rightarrow x$, $\zeta_n\stackrel{w}\longrightarrow \zeta$ and $\zeta_n\in N^p(S, x_n)$ for each $n\in \mathbb{N}$.

Let $f:H\rightarrow \mathbb{R}\cup\{+\infty\}$ be a lower semicontinuous function. We denote
$$
{\rm dom}(f):=\{y\in H: f(y)<+\infty\}\,\,\,{\rm and}\,\,\,{\rm epi}f:=\{(x, \alpha)\in X\times \mathbb{R}: f(x)\leq \alpha\}
$$
the {\it domain} and the {\it epigraph} of $f$, respectively. Let $x\in {\rm dom}(f)$. Recall that $f$ is said to be {\it G\^ateaux differentiable} at $x$, provided there exists $\xi\in H$ such that $$
\lim_{t\downarrow 0}\frac{f(x+th)-f(x)}{t}=\langle \xi, h\rangle,\,\,\,\forall\,h\in H.
$$
Recall from \cite{8} that a vector $\zeta\in H$ is said to be a  {\em proximal
subgradient} of $f$ at $x$ provided that
\begin{equation*}
(\zeta, -1)\in N^p({\rm epi}f, (x, f(x))).
\end{equation*}
The set of all such $\zeta$, denoted $\partial_pf(x)$, is said to be the {\em proximal subdifferential} of $f$ at $x$. It is known from \cite{8} that one useful and important characterization for proximal subgradient is as follows: $\zeta\in \partial_pf(x)$ if and only if there exist two constants $r, \delta\in ]0, +\infty[$ such that
\begin{equation*}
f(y)\geq f(x)+\langle\zeta,y-x\rangle-\frac{r}{2}\|y-x\|^2,\,\,\,\forall~
y\in B(x,\delta).
\end{equation*}

The {\it limiting(Mordukhovich) subdifferential} of $f$ at $x$ is defined as:
\begin{equation*}
\partial f(x):=\{\zeta\in H: (\zeta, -1)\in N({\rm epi}f, (x, f(x)))
\}.
\end{equation*}
It is proved in \cite{14,15} that
$$
\partial f(x)=\mathop{\rm Limsup}_{y\stackrel{f}\rightarrow x}\partial_pf(y),
$$
where $y\stackrel{f}\rightarrow x$
means $y\rightarrow x$ and $f(y)\rightarrow f(x)$; that is, $\zeta\in\partial f(x)$ if and only if there exist $x_n\stackrel{f}\rightarrow x$ and $\zeta_n\stackrel{w}\longrightarrow\zeta$ such that $\zeta_n\in\partial_pf(x_n)$ for all $n\in\mathbb{N}$. When $f$ is convex, the proximal subdifferential and the limiting subdifferential of $f$ at $x$ coincide and both reduce the subdifferential in the sense of convex analysis, that is
$$
\partial_pf(x)=\partial f(x)=\{\zeta\in H: \langle \zeta, y-x\rangle\leq f(y)-f(x), \ \ \forall y\in H\}.
$$
Readers are invited to consult books \cite{8,14,20} for more details on these various normal cones and subdifferentials.

In the analysis, we need the notion of {\it upper semicontinuity}. Let $X, Y$ be Hausdorff topological spaces and let $F:X\rightrightarrows Y$ be a set-valued mapping from $X$ into $Y$. Recall from \cite{AF} that $F$ is {\it upper semicontinuous} at $x\in X$ iff, for each open set $V$ in $Y$ containing $F(x)$, there is an open neighborhood $U$ of $x$ such that $F(x')\subset V$ for all $x'\in U$.\\

Let $f:\mathbb{R}^m\rightarrow \mathbb{R}$ be a smooth function. We consider the following unconstraint optimization problem:
\begin{equation}\label{1}
\min_{x\in \mathbb{R}^m} f(x).
\end{equation}
Note that it is of significance to solve it by using the iteration of line search methods. Each iteration of a line search method computes a search direction $u_{n-1}$ and then decides how far to move along this direction. The iteration is given by
$$
x_n:=x_{n-1}+\varepsilon_{n-1}u_{n-1},
$$
where the positive scalar $\varepsilon_{n-1}$ is called the step length. When it comes to problem \eqref{1}, it is used generally to find the steepest descent unit direction, and then $f$ can be reduced along this direction and finally attain the global minimizer hopefully. The classical algorithm for smooth functions by the steepest descent method is as follows:

{\it {\bf{Algorithm SDM}} Let $f:\mathbb{R}^m\rightarrow \mathbb{R}$ be smooth and let $x_0$ be given.

{\rm Step 1}: Let $n:=1$. Go to {\rm Step 2}.

{\rm Step 2}: Check vector $\triangledown f(x_{n-1})$.  If $\triangledown f(x_{n-1})=0$, then terminate; otherwise, set
\begin{equation*}
u_{n-1}:=-\frac{\triangledown f(x_{n-1})}{\|\triangledown f(x_{n-1})\|}\ \ and\ \ x_n:=x_{n-1}+\varepsilon_{n-1}u_{n-1},
\end{equation*}
where the step length $\varepsilon_{n-1}$ is chosen by
$$
f(x_{n-1}+\varepsilon_{n-1}u_{n-1})=\min_{t\geq 0} f(x_{n-1}+tu_{n-1}).
$$
Go to {\rm Step 3}.

{\rm Step 3}: Set $n:=n+1$ and return to {\rm Step 2}.}\\

For the smooth function case, the sequence $u_{n}$, produced by Algorithm SDM, can converge to a local minimizer only by making additional requirements on the search direction $u_{n}$. However, the following known convergent result shows that, for twice differentiable convex functions, the sequence generated by Algorithm SDM would converge to the global minimizer with no assumption on search direction (cf. \cite[Theorem 2.3.8]{P}):\\

\noindent{\bf Theorem 2.1.} {\it Let $f:\mathbb{R}^m\rightarrow \mathbb{R}$ be a strictly convex, twice continuously differentiable function and $x_0\in\mathbb{R}^m$. Denote $\Omega:=\{x\in\mathbb{R}^m: f(x)\leq f(x_0)\}$ the level set. Suppose that there exist constants $m, M\in ]0, +\infty[$ such that
\begin{equation}\label{2}
  m\|y\|^2\leq\langle y, \triangledown^2f(z)y\rangle\leq M\|y\|^2\ \ \forall z\in\Omega\ {\rm and}\ \forall y\in \mathbb{R}^m.
\end{equation}
Let $\{x_n\}$ be generated by Algorithm SDM. Then there exists a subsequence  $\{x_{n_k}\}$ of $\{x_n\}$ converging to a point $\bar x$ that minimizes $f(z)$ over $z\in\mathbb{R}^m$, where $\triangledown^2f(z)$ denotes the Hessian matrix of $f$ at $z$.
}

Note that qualification condition \eqref{2} could guarantee the boundedness of the level set $\Omega$ (by \cite[Theorem B2.8]{P}) and the bounds $m, M$ may be the smallest and largest eigenvalue of the Hessian matrix $\triangledown^2f(z)$. Our main work in this paper is to establish the corresponding algorithm for nonsmooth functions defined on Hilbert space and use it to find stationary points of those functions satisfying prox-regularity and Lipschitzian continuity (see Theorem 3.1). Then, we apply the newly-established algorithm to the search of minimizer for lower semicontinuous G\^{a}teaux differentiable convex functions defined on the finite-dimensional space so as to extend and improve Theorem 2.1 above.

\section{Nonsmooth Steepest Descent Method}
In this section, we study nonsmooth steepest descent method for nonsmooth functions on the Hilbert space and establish the corresponding algorithm by proximal subgradients so as to extend the original smooth algorithm by steepest descent method and steepest descent gradient. Then, we use the newly-established algorithm to find stationary points for those nonsmooth functions satisfying prox-regularity and Lipschitzian continuity. We begin with recalling the definition of {\it prox-regularity} taken from \cite{20,19}.\\

\noindent{\bf Definition 3.1.} Let $f:H\rightarrow \mathbb{R}\cup\{+\infty\}$, and let $\bar
y\in {\rm dom}(f)$. We say that $f$ is {\it prox-regular} at $\bar y$ for
$\bar{w}\in \partial_pf(\bar y)$ iff there exist $\delta>0$ and
$L\geq 0$ such that
\begin{equation}\label{3}
f(y')\geq f(y)+\langle w,
y'-y\rangle-\frac{L}{2}\|y'-y\|^2\,\,\,\,\forall y'\in B(\bar y,
\delta),
\end{equation}
whenever $y\in B(\bar y, \delta)$ with $|f(y)-f(\bar y)|<\delta$ and
$w\in \partial_pf(y)\cap B(\bar{w}, \delta)$. We say that $f$ is
prox-regular at $\bar y$ iff this holds for
every $\bar w\in \partial_pf(\bar y)$.\\

The following concept introduced in \cite{5} is also relevant here.\\

\noindent {\bf Definition 3.2.} Let $f:H\rightarrow \mathbb{R}\cup\{+\infty\}$, and let $E\subset
H$. We say that $f$ is {\it uniformly prox-regular} on $E$ iff there exist
$\delta>0$ and $L>0$ such that, for any $\bar y\in E$ and $\bar w\in
\partial_pf(\bar y)$, one has
\begin{equation}\label{4}
f(y')\geq f(y)+\langle w,
y'-y\rangle-\frac{L}{2}\|y'-y\|^2\,\,\,\,\forall y'\in B(\bar y,
\delta),
\end{equation}
whenever $y\in B(\bar y, \delta)$ with $|f(y)-f(\bar y)|<\delta$ and
$w\in \partial_pf(y)\cap B(\bar{w}, \delta)$. We say that $f$ is
{\it locally uniformly prox-regular} around $y_0\in H$ iff $E$ can be taken as a neighborhood of $y_0$, i.e., $E=B(y_0, r)$ for some $r>0$.\\

The following proposition, obtained from \cite[Proposition 3.3]{5}, is one characterization for the local uniform prox-regularity.
\begin{pro}
Let $f:H\rightarrow \mathbb{R}\cup\{+\infty\}$ and $y_0\in
H$. Then $f$ is locally uniformly prox-regular around $y_0$ if and
only if there are $\delta>0$ and $L>0$ such that, for any $y,
y'\in B(y_0, \delta)$ and $w\in \partial_pf(y)$, we have the estimate
\begin{equation}\label{5}
f(y')\geq f(y)+\langle w, y'-y\rangle-\frac{L}{2}\|y'-y\|^2.
\end{equation}
\end{pro}

Now, we provide one algorithm of nonsmooth steepest descent method for nonsmooth functions by proximal subgradients and use it to find stationary points for nonsmooth functions with prox-regularity and Lipschitzian continuity assumptions in the Hilbert space.
\begin{them}
Let $f:H\rightarrow \mathbb{R}\cup\{+\infty\}$ be lower
semicontinuous and bounded below. Let $x_0\in {\rm{dom}}(f)$ with
$\partial_p f(x_0)\neq\emptyset$. We define the level set as
\begin{equation}\label{6}
\Omega:=\{x\in H : f(x)\leq f(x_0)\}.
\end{equation}
Suppose that there exist $\delta, L>0$ such that $f$ is uniformly prox-regular on the set $B(\Omega, \delta):=\{x\in H: d(x, \Omega)<\delta\}$ and
\begin{equation}\label{7}
\partial_p f(y)\subset\partial_p f(z)+L\|y-z\|B_H\,\,\,\,\,\,\forall\, y\in\Omega\,\,{\rm with}\,\,\,\,\|y-z\|<\delta.
\end{equation}
We use the following Algorithm NSDM to construct a sequence $\{x_n:n=1,2,\cdots\}$:

{\bf{Algorithm NSDM}} Let $x_0$ be given with $\partial_p f(x_0)\neq\emptyset$.

{\rm Step 1}: Let $n:=1$. Go to {\rm Step 2}.

{\rm Step 2}: Choose $v_{n-1}\in \partial_p f(x_{n-1})$ and check vector $v_{n-1}$.  If $v_{n-1}=0$, then terminate; otherwise, set
$$
u_{n-1}:=-\frac{v_{n-1}}{\|v_{n-1}\|}\,\,\,\,{\rm and}\,\,\,\,x_n:=x_{n-1}+\varepsilon_{n-1}u_{n-1},
$$
where $\varepsilon_{n-1}$ is chosen as follows:
\begin{equation}\label{8}
f(x_{n-1}+\varepsilon_{n-1}u_{n-1})=\min_{t\geq 0}f(x_{n-1}+tu_{n-1}).
\end{equation}
Go to {\rm Step 3}.

{\rm Step 3}: Set $n:=n+1$ and return to {\rm Step 2}.

Then Algorithm NSDM is valid and every cluster point of $\{x_n:n=0,1,\cdots\}$ is the stationary point.
\end{them}

\noindent{\it Proof} First, we have to show that $\Omega\subset{\rm dom}(\partial_p f)$ so that the validness of Algorithm NSDM could be  guaranteed. Let $y\in \Omega$. By the Density Theorem (cf. \cite[Theorem 3.1]{8}), one has ${\rm dom}(\partial_p f)\cap{\Omega}$ is dense in $\Omega$. Then, there exist $y_k\in{\rm dom}(\partial_p f)\cap\Omega$ such that $y_k\rightarrow y$. By virtue of \eqref{7}, when $k$ is sufficiently large, one has
$$
\partial_p f(y_k)\subset\partial_p f(y)+L\|y_k-y\|B_H.
$$
This implies that $\partial_p f(y)\neq\emptyset$ as $\partial_p f(y_k)\neq\emptyset$ for all $k$ sufficiently large.
 Hence $\Omega\subset{\rm dom}(\partial_p f)$.

Since $f$ is uniformly prox-regular, by Proposition 3.1, without any
loss of generality, we can assume that, for any $\bar y\in B(\Omega, \delta)$
and $\bar{w}\in \partial_pf(\bar y)$, one has
\begin{equation}\label{9}
f(y')\geq f(y)+\langle w, y'-y\rangle-\frac{L}{2}\|y'-y\|^2\,\,\,\,\forall
y'\in B(\bar y, \delta)
\end{equation}
whenever $w\in \partial_pf(y)\cap B(\bar{w}, \delta)$ and $y\in
B(\bar y, \delta)$ with $|f(y)-f(\bar y)|<\delta$, perhaps after
reducing the value of $\delta>0$ and increasing the value of $L>0$.
If there exists $n_0$ such that $v_{n_0}=0$, then the conclusion
holds. Next, we suppose that $v_{n}\neq 0$ for
all $n\in \mathbb{N}$. Let $\varepsilon\in ]0, \delta[$. Noting that $v_{n-1}\in
\partial_p f(x_{n-1})$, it follows from \eqref{7} that there exist
$w_{n-1}\in \partial_pf(x_{n-1}+\varepsilon u_{n-1})$ and $b_{n-1}\in B_H$ such
that
\begin{equation}\label{10}
v_{n-1}=w_{n-1}+L\varepsilon b_{n-1}
\end{equation}
(thanks to $\varepsilon\in ]0, \delta[$). By \eqref{9} and \eqref{10}, one has
$$
f(x_{n-1})\geq f(x_{n-1}+\varepsilon u_{n-1})+\langle w_{n-1}, -\varepsilon u_{n-1}\rangle-\frac{L}{2}\varepsilon^2
$$
$$
\,\,\,\,\,\,\,\,\,\,\,\,\,\,\,\,\,\,\,\,\,\,\,\,\,\,\,\,\,\,\,\,\,\,\,\,\geq f(x_{n-1}+\varepsilon u_{n-1})-\varepsilon\langle v_{n-1},  u_{n-1}\rangle-L\varepsilon^2-\frac{L}{2}\varepsilon^2.
$$
This and \eqref{8} imply that
\begin{equation}\label{11}
\|v_{n-1}\|\leq \frac{f(x_{n-1})-f(x_n)}{\varepsilon}+\frac{3}{2}L\varepsilon.
\end{equation}
Since $\{f(x_n)\}$ is monotonic decreasing by \eqref{8} and $f$ is bounded below, it follows that $\lim_{n\rightarrow \infty}f(x_n)$ exists. Hence, there exists $N\in \mathbb{N}$ such that when $n>N$, one has
$$
\frac{f(x_{n-1})-f(x_n)}{\varepsilon}<\varepsilon.
$$
From this and \eqref{11}, we have
$$
\|v_{n-1}\|\leq(1+\frac{3}{2}L)\varepsilon,\,\,\,\,\,n>N,
$$
and thus $v_n\rightarrow 0$. Let $\bar x$ be any cluster point of $\{x_n:n=0,1,\cdots\}$. Then there exists a subsequence $\{x_{n_k}\}$ of $\{x_n\}$ such that $x_{n_k}\rightarrow\bar x \in H$. It follows from the lower semicontinuity of $f$ that
$$
f(\bar x)\leq\liminf_{k\rightarrow\infty}f(x_{n_k})\leq f(x_0),
$$
and hence $\bar x\in \Omega$. By applying \eqref{9} again, for $k$ large enough, one has
$$
f(\bar x)\geq f(x_{n_k})+\langle v_{n_k}, \bar x-x_{n_k}\rangle-\frac{L}{2}\|\bar x-x_{n_k}\|^2,
$$
and consequently $f(\bar x)\geq\limsup_{k\rightarrow \infty}f(x_{n_k})$. This and lower semicontinuity of $f$ imply that $f(x_{n_k})\rightarrow f(\bar x)$ as $k\rightarrow \infty$. Noting that $v_{n_k}\in \partial_pf(x_{n_k})$ and $v_{n_k}\rightarrow 0$, it follows that $0\in \partial f(\bar x)$. The proof is complete. \hfill $\Box$\\

\noindent{\bf Remark 3.1.} (i) Condition \eqref{7} in Theorem 3.1 refers to the Lipschitzian property of subdifferential $\partial_pf$. Recently  Bacak, Borwein, Eberhard and Mordukhovich \cite{B} mainly studied the Lipschitzian property of subdifferentials for prox-regular functions in Hilbert space, and obtained results on this Lipschitzian property of subdifferentials (cf. \cite[Lemma 5.2 and Theorem 5.3]{B}).  For the smooth functions, subdifferential $\partial_pf$ reduces to Fr\'{e}chet derivative operator $\triangledown f$, and condition \eqref{7} means that  $\triangledown f$ is Lipschitzian continuous with modulus $L>0$. This condition holds for $C^{1, 1}$ functions, para-convex and para-concave functions and twice differentiable functions on finite-dimensional space with bounded Hessian matrices. In particular, the convex function studied in Theorem 2.1 satisfies condition \eqref{7}.

(ii) Algorithm NSDM provides a search direction $u_{n-1}$ by proximal subgradient for nonsmooth functions defined on the Hilbert space and generates a sequence $\{x_n\}$ having some subsequence converging the stationary point with mild assumptions. When restricted to the smooth case, Algorithm NSDM reduces to Algorithm SDM mentioned in Section 2. \\

The following example shows that Algorithm NSDM is practicable to find stationary points (or even minimizers) of functions satisfying prox-regularity and Lipschitzian continuity assumptions in Hilbert space.\\

\noindent{\bf Example 3.1.} Let $H:=l^2$ and $y=(\eta(1),\cdots,\eta(k),\cdots)\in H$. Define a function $f$ on $H$ by
$$
f(x):=\frac{1}{2}\sum_{i=1}^{\infty}\xi(i)^2-\sum_{i=1}^{\infty} \xi(i)\eta(i)+3, \ \ {\rm for \ all \ }x=(\xi(1),\cdots,\xi(k),\cdots)\in H.
$$
Clearly $f$ is convex, and thus satisfies the uniform prox-regularity on $H$. Since $\partial_pf(x)=x-y$ (by computing), then $\partial_pf$ is Lipschitzian continuous on $H$. Take $x_0:=(0,\cdots,0,\cdots)\in H$ and let $\Omega$ be defined as \eqref{6}. One can verify that \eqref{7} holds with $L:=1$. Using Algorithm NSDM, the sequence $\{x_n: n=1,2,\cdots\}$ can be generated and $v_n=x_n-y$. By virtue of the proof of Theorem 3.1, sequence $\{v_n\}$ converges to vector $0\in H$, and thus $y$ is a cluster point of $\{x_n\}$. It follows from Theorem 3.1 that $y$ is a stationary point, and consequently $y$ minimizes $f(x)$ over $x\in H$ due to the convexity of $f$.\\

Let $F: H\rightrightarrows H$ be a set-valued mapping and $x\in H$ with $F(x)\not=\emptyset$. Recall from \cite{AF,14} that $F$ is said to be {\it locally Lipschitzian} at $x$ iff there exist $\delta, \kappa>0$ such that
$$
F(x')\subset F(x'')+\kappa\|x'-x''\|B_H\,\,\,\,\,\forall\,x', x''\in B(x, \delta).
$$
\begin{them}
Let $f:\mathbb{R}^m\rightarrow \mathbb{R}\cup\{+\infty\}$ be
lower semicontinuous and bounded below. Let $x_0\in {\rm{dom}}(f)$
with $\partial_p f(x_0)\neq\emptyset$ and the level set $\Omega$ be
defined as \eqref{6}. Suppose that $\Omega$ is bounded, $\partial_pf$ is locally Lipschitzian on $\Omega$ and that $f$ is locally uniformly
prox-regular on $\Omega$. Let the sequence $\{x_n: n=0,1,\cdots\}$
be generated by Algorithm NSDM. Then there exists a subsequence $\{x_{n_k}\}$ of $\{x_n\}$ such that $\{x_{n_k}\}$ converges to a point $\bar x\in \Omega$ with $0\in \partial f(\bar x)$.
\end{them}

\noindent{\it Proof} First, we need to show that $\Omega\subset{\rm dom}(\partial_p f)$
so that the Algorithm NSDM is valid. Let $y\in \Omega$. By the
Density Theorem, one has ${\rm dom}(\partial_p f)\cap{\Omega}$ is
dense in $\Omega$. Then, there exists $y_k\in{\rm dom}(\partial_p
f)\cap\Omega$ such that $y_k\rightarrow y$. Noting that
$\partial_pf$ is locally Lipschitzian at $y$, it follows that there exist $r>0$
and $\kappa_y>0$ such that
$$
\partial_pf(z_1)\subset\partial_pf(z_2)+\kappa_y\|z_1-z_2\|B_H\,\,\,\,\forall
z_1, z_2\in B(y, r).
$$
This means that when $k$ is sufficiently large, one has
$$
\partial_p f(y_k)\subset\partial_p f(y)+\kappa_y\|y_k-y\|B_H,
$$
and consequently $\partial_p f(y)\neq\emptyset$ as $\partial_p
f(y_k)\neq\emptyset$ for all $k$ sufficiently large. Hence $\Omega\subset{\rm
dom}(\partial_p f)$.

Let $y\in \Omega$. Since $\partial_pf$ is locally Lipschitzian
and locally uniformly prox-regular around $y$, by Proposition 3.1,
there exist $\delta_y>0$ and
$L_y>0$ such that, for any $y', y''\in B(y, \delta_y)$ and $w\in
\partial_pf(y')$, one has
\begin{equation}\label{12}
f(y'')\geq f(y')+\langle w, y''-y'\rangle-\frac{L_y}{2}\|y''-y'\|^2,
\end{equation}
and
\begin{equation}\label{13}
\partial_pf(z_1)\subset\partial_pf(z_2)+L_y\|z_1-z_2\|B_H\,\,\,\,\forall
z_1, z_2\in B(y, \delta_y),
\end{equation}
perhaps after reducing the value of $\delta_y>0$ and increasing the
value of $L_y>0$ if necessary. Then
\begin{equation}\label{14}
\Omega\subset\bigcup\limits_{y\in \Omega}B(y, \frac{\delta_y}{2}).
\end{equation}
Note that $\Omega$ is closed by the semicontinuity of $f$ and
bounded, so $\Omega$ is compact. It follows from \eqref{14}
that there exist $y_1,\cdots,y_l\in\Omega$ such that
\begin{equation}\label{15}
\Omega\subset\bigcup_{i=1}^lB(y_i, \frac{\delta_{y_i}}{2}).
\end{equation}
Let
$$\delta:=\min\Big\{\frac{\delta_{y_1}}{2},\cdots,\frac{\delta_{y_l}}{2}\Big\}\,\,\,\,{\rm
and}\,\,\,\,L:=\max\{L_{y_1},\cdots,L_{y_l}\},
$$
and let
$\varepsilon\in ]0, \delta[$. By the Algorithm NSDM, we have
$x_{n-1}\in\Omega$, and thus there exists $i\in\{1, \cdots, l\}$ such that
$x_{n-1}\in B(y_{i}, \frac{\delta_{y_i}}{2})$. Noting that
$v_{n-1}\in
\partial_pf(x_{n-1})$ and $x_{n-1}+\varepsilon u_{n-1}\in B(y_{i},
\delta_{y_i})$, it follows from \eqref{13} that
$$
v_{n-1}\in\partial_pf(x_{n-1})\subset\partial_pf(x_{n-1}+\varepsilon
u_{n-1})+L\varepsilon B_{\mathbb{R}^m}.
$$
(thanks to $L\geq L_{y_i}$). Then, there are $w_{n-1}\in
\partial_pf(x_{n-1}+\varepsilon u_{n-1})$ and $b_{n-1}\in B_{\mathbb{R}^m}$ such that
\begin{equation}\label{16}
v_{n-1}=w_{n-1}+L\varepsilon b_{n-1}
\end{equation}
Since $x_{n-1}, x_{n-1}+\varepsilon u_{n-1}\in B(y_i, \delta_{y_i})$
and $w_{n-1}\in
\partial_pf(x_{n-1}+\varepsilon u_{n-1})$, by \eqref{12}
and \eqref{16}, one has
$$
f(x_{n-1})\geq f(x_{n-1}+\varepsilon u_{n-1})+\langle w_{n-1},
-\varepsilon u_{n-1}\rangle-\frac{L_{y_i}}{2}\varepsilon^2
$$
$$
\,\,\,\,\,\,\,\,\,\,\,\,\,\,\,\,\,\,\,\,\,\,\,\,\,\,\,\,\geq
f(x_{n-1}+\varepsilon u_{n-1})-\varepsilon\langle v_{n-1},
u_{n-1}\rangle-L\varepsilon^2-\frac{L}{2}\varepsilon^2.
$$
This and \eqref{8} imply that
$$
\|v_{n-1}\|\leq\frac{f(x_{n-1})-f(x_n)}{\varepsilon}+\frac{3}{2}L\varepsilon.
$$
Since $\{f(x_n)\}$ is monotonic decreasing and bounded below, one
has $\lim_{n\rightarrow \infty}f(x_n)$ exists, and
consequently there exists $N\in \mathbb{N}$ such that, whenever
$n>N$, one has
$$
\frac{f(x_{n-1})-f(x_n)}{\varepsilon}<\varepsilon.
$$
Hence
$$
\|v_{n-1}\|<(1+\frac{3}{2}L)\varepsilon,\,\,\,\,\,\forall n>N,
$$
and consequently $v_n\rightarrow 0$. Noting that $\{x_n\}\subset
\Omega$ and $\Omega$ is compact, there exists a subsequence $\{x_{n_k}\}$ of  $\{x_n\}$ converging to a point $\bar x\in \Omega$. By using \eqref{12}, when $k$ is sufficiently large, we have
$$
f(\bar x)\geq f(x_{n_k})+\langle v_{n_k}, \bar x-x_{n_k}\rangle-\frac{L}{2}\|\bar x-x_{n_k}\|^2.
$$
Thus $f(\bar x)\geq \limsup_{k\rightarrow \infty}f(x_{n_k})$, and consequently $f(x_{n_k})\rightarrow f(\bar x)$ as $f$ is lower semicontinuous at $\bar x$. Hence $0\in \partial
f(\bar x)$. The proof is complete. \hfill$\Box$\\

Next, we use the established Algorithm NSDM to find the minimizer of lower semicontinuous and convex functions defined on finite-dimensional space. Since convexity implies local uniform prox-regularity automatically, the following theorem is immediate from Theorem 3.2.

\begin{them}
Let $f:\mathbb{R}^m\rightarrow \mathbb{R}\cup\{+\infty\}$ be a
lower semicontinuous, convex and bounded below function, $x_0\in {\rm{dom}}(f)$ with $\partial f(x_0)\neq\emptyset$
and the level set $\Omega$ be
defined as \eqref{6}. Suppose that $\Omega$ is bounded and $\partial f$ is locally Lipschitzian on $\Omega$. Let the sequence $\{x_n: n=0,1,\cdots\}$
be generated by Algorithm NSDM. Then there exists a subsequence $\{x_{n_k}\}$ of $\{x_n\}$ such that $\{x_{n_k}\}$ converges to a point $\bar x\in \Omega$ with $0\in \partial f(\bar x)$, and consequently $f(\bar x)=\min_{z\in \mathbb{R}^m}f(z)$.
\end{them}

Now, we provide one convergent theorem for lower semicontinuous, convex and G\^{a}teaux differentiable functions defined on the finite-dimensional space.  This theorem shows that the sequence generated by Algorithm NSDM converges to the minimizer even without Lipschitzian assumption of subdifferential.

\begin{them}
Let $f:\mathbb{R}^m\rightarrow \mathbb{R}\cup\{+\infty\}$ be lower semicontinuous, convex, G\^{a}teaux differentiable and bounded below function, $x_0\in {\rm{dom}}(f)$
and the level set $\Omega$ be
defined as \eqref{6}. Suppose that $\Omega$ is bounded and $\Omega\subset{\rm int}({\rm dom}(f))$. Let $\{x_n\}$ be generated by the Algorithm NSDM. Then there exists a subsequence $\{x_{n_k}\}$ of $\{x_n\}$ such that $\{x_{n_k}\}$ converges to a point $\bar x\in \Omega$ that minimizes $f(z)$ over $z\in \mathbb{R}^m$.
\end{them}

\noindent{\it Proof}
First, the validness of the Algorithm NSDM follows from the G\^{a}teaux differentiability of $f$.

Using \cite[Proposition 3.3]{18}, we have that $f$ is continuous on ${\rm int}({\rm dom}(f))$, and consequently the subdifferential mapping $\partial f$ is norm-to-norm upper semicontinuous on ${\rm int}({\rm dom}(f))$ by \cite[Proposition 2.5]{18}. Let $\varepsilon\in ]0, +\infty[$. Noting that $\mathbb{R}^m$ is finite-dimensional and $\Omega\subset{\rm int}({\rm dom}(f))$, by the upper semicontinuity, for each $y\in\Omega$ there exists $\delta_y>0$ such that
\begin{equation}\label{17}
\partial f(y')\subset\partial f(y)+\frac{\varepsilon}{3} B_{\mathbb{R}^m}\,\,\,\,\forall y'\in B(y, \delta_y).
\end{equation}
Note that $\Omega$ is bounded, and consequently $\Omega$ is compact. This implies that there exist $y_1,\cdots, y_l\in\Omega$ such that
\begin{equation}\label{18}
\Omega\subset\bigcup_{i=1}^lB(y_i, \frac{\delta_{y_i}}{2}).
\end{equation}
Let $\delta:=\min\Big\{\frac{\delta_{y_1}}{2},\cdots,\frac{\delta_{y_l}}{2}\Big\}$ and $\alpha\in ]0, \delta[$. By the Algorithm NSDM, we have
$x_{n-1}\in\Omega$, and there exists $i\in\{1,\cdots,l\}$ such that
$x_{n-1}\in B(y_{i}, \frac{\delta_{y_i}}{2})$. Since $f$ is G\^{a}teaux differentiable at $y_i$, then the subdifferential $\partial f(y_i)$ is the singleton. Noting that $v_{n-1}\in \partial f(x_{n-1})$ and $x_{n-1}+\alpha u_{n-1}\in B(y_i, \delta_{y_i})$, it follows from \eqref{18} that there exist $w_{n-1}\in\partial f(x_{n-1}+\alpha u_{n-1})$ and $a_{n-1}, b_{n-1}\in B_{\mathbb{R}^m}$ such that
$$
v_{n-1}=w_{n-1}+\frac{\varepsilon}{3}(b_{n-1}+a_{n-1}).
$$
This and convex inequality in convex analysis imply that
$$
f(x_{n-1})-f(x_{n-1}+\alpha u_{n-1})\geq\langle w_{n-1}, -\alpha u_{n-1} \rangle\geq \langle v_{n-1}, -\alpha u_{n-1}\rangle-\frac{2}{3}\alpha\varepsilon
$$
and, by \eqref{8}, one has
\begin{equation}\label{19}
\|v_{n-1}\|\leq\frac{f(x_{n-1})-f(x_n)}{\alpha}+\frac{2}{3}\varepsilon.
\end{equation}
Since $\{f(x_n)\}$ is monotonic decreasing and bounded below, then $\lim_{n\rightarrow \infty}f(x_n)$ exists and there exists $N\in \mathbb{N}$ such that, whenever $n>N$, we have
$$
\frac{f(x_{n-1})-f(x_n)}{\alpha}<\frac{\varepsilon}{3}.
$$
From this and \eqref{19}, one has
$$
\|v_{n-1}\|<\varepsilon,\ \ \forall n>N.
$$
This implies that $v_{n}\rightarrow 0$. Since $\{x_n\}\subset \Omega$ and $\Omega$ is compact, there exists a subsequence $\{x_{n_k}\}$ of $\{x_n\}$ such that $x_{n_k}\rightarrow \bar x\in \Omega$. Note that $f(x_{n_k})\rightarrow f(\bar x)$ follows from the continuity of $f$ at $\bar x$, so $0\in \partial f(\bar x)$ and this is equivalent to that $f(\bar x)=\min_{z\in \mathbb{R}^m}f(z)$ (thanks to the convexity of $f$). The proof is complete. \hfill$\Box$\\

\noindent{\bf Remark 3.2.} In Theorem 3.4, Algorithm NSDM established has been used to construct a sequence having some subsequence converging to the minimizer of lower semicontinuous convex and G\^{a}teaux differentiable functions that are more general than strictly convex and twice continuously differentiable functions. Hence this convergent result improves Theorem 2.1 in the sense of relaxing the twice continuously differentiable assumption.

\section{Concluding Remarks}
This paper is devoted to presenting an algorithm for finding stationary points of nonsmooth functions based on nonsmooth steepest descent method and proximal subgradients in Hilbert space. The algorithm established is used to produce one sequence converging to zero vector so as to obtain the stationary point of lower semicontinuous functions satisfying prox-regularity and Lipschitzian continuity on Hilbert space. This algorithm is also applied to the search of minimizers for convex functions on finite-dimension space. The obtained convergent theorem extends the existing result on twice continuously differentiable convex functions on finite-dimensional space.\\

%The authors are indebted to two anonymous referees for their helpful comments which allowed us to improve the original presentation.
\noindent{\bf Acknowledgement} This research was supported by the National Natural Science Foundation of P. R. China (Grant No. 11261067), the Scientific Research Foundation of Yunnan University under grant No. 2011YB29 and by IRTSTYN.

%{\bf Acknowledgment.} The authors wish to thank the referees for many valuable comments and for
%reference [1,4,9] as well as drawing authors' attention to amenable functions due to Poliquin and
%Rockafellar and penalty functions.


\begin{thebibliography}{99}


\bibitem{3} Attouch, H., Cominetti, R.: A dynamical approach to convex minimization coupling approximation with the steepest descent method, J. Differ. Equations., 128, 519-540 (1996).

\bibitem{4} Bazaraa, M. S., Shetty, C. M.: Nonlinear programming: theory and algorithm, New York: Wiley (1979).

\bibitem{CY} Chuong, T. D., Yao, J.-C.: Steepest descent methods for critical points in vector optimization problems, Optimization, 91, 1811-1829 (2012).

\bibitem{DS} Grana Drummond, L. M., Svaiter, B. F.: A steepest descent method for vector optimization, J. Comput. Appl. Math., 175, 395-414 (2005).

\bibitem{FS} Fliege, J., Svaiter, B. F.: Steepest descent methods for multicriteria optimization, Math. Methods Oper. Res., 51, 479-494 (2000).

\bibitem{17} Nocedal, J., Wright, S. J.: Numerical optimization, Springer, New York (2006).

\bibitem{P} Polak, E.: Computational methods in optimization, Academic Press, New York (1971).

\bibitem{23} Smyrlis, G., Zisis, V.: Local convergence of the steepest descent method in Hilbert spaces, J. Math. Anal. Appl., 300, 436-453 (2004).

\bibitem{10} Hestenes, M. R., Stiefel, E.: Methods of conjugate gradients for solving linear systems, Journal of Research of the National Bureau of standards, 49, 409-436 (1952).

\bibitem{9} Fletcher, R., Reeves, C. M.: Functions minimization by conjugate gradients, Comput. J., 7, 149-154 (1964).


\bibitem{8} Clarke, F. H., Ledyaev, Yu. S., Stern, R. J., Wolenski, P. R.: Nonsmooth analysis and control theory, Springer (1998).

\bibitem{14} Mordukhovich, B. S.: Variational Analysis and Generalized differentiation I/II, Springer-verlag, Berlin (2006).

\bibitem{15} Mordukhovich, B. S., Shao, Y.: Nonsmooth sequential analysis in Asplund spaces, Trans. Amer. Math. Soc., 348, 1235-1280 (1996).

\bibitem{20} Rockafellar, R. T., Wets, R. J-B.: Variational Analysis, Springer-Verlag, New York (1998).





\bibitem{AF} Aubin, J. P., Frankowska, H.: Set-valued Analysis, Birkh\"{a}user, Boston (1990).

\bibitem{19} Poliquin, R., Rockafellar, R.T.: Prox-regular functions in variational analysis, Trans. Amer. Math. Soc., 348, no. 5, 1805-1838 (1996).



\bibitem{5} Bernard, F., Thilbaut, L.: Uniform prox-regularity of functions and epigraphs in Hilbert spaces, Nonlinear Anal., 60, 187-207 (2005).

\bibitem{B} Bacak, M., Borwein, J. M., Eberhard, A., Mordukhovich, B. S.: Infimal convolutions and Lipschitzian properties of subdifferentials for prox-regular functions in Hilbert spaces,  J. Convex Anal., 17, 737-763 (2010).

\bibitem{18} Phelps, R. R.: Convex functions, Monotone operators and Differentiability, Lecture Notes in Math. 1364, Springer, New York (1989).

















\end{thebibliography}
\end{document}